\documentclass[12pt,reqno]{amsart} 

\usepackage{amssymb,latexsym}
\usepackage[draft=true]{hyperref}
\usepackage{cite} 

\usepackage[height=190mm,width=130mm]{geometry} 

\theoremstyle{plain}
\newtheorem{theorem}{Theorem}
\newtheorem{lemma}{Lemma}

\theoremstyle{definition}

\theoremstyle{remark}

\numberwithin{equation}{section} 
\begin{document}
\title[Para-Sasakian Manifolds and $*$-Ricci Solitons]{Para-Sasakian Manifolds and $*$-Ricci Solitons} 

\author{D. G. Prakasha}
\address{Department of Mathematics\\ Karnatak University\\ Dharwad - 580003\\India}

\email{prakashadg@email.com, dgprakasha@kud.ac.in, viru0913@gmail.com}

\thanks{The first author (DGP) was supported in part by University Grants Commission, New Delhi, India in the form of UGC-SAP-DRS-III programme to the Department of Mathematics, K.U. Dharwad.}

\author{Pundikala Veeresha}

\begin{abstract}
In this paper we study a special type of metric
called $*$-Ricci soliton on para-Sasakian manifold. We prove that if the para-Sasakian metric is a $*$-Ricci soliton on a manifold $M$, then $M$ is either $\mathbb{D}$-homothetic to an Einstein manifold, or  the Ricci tensor of $M$ with respect to the canonical paracontact connection vanishes.
\end{abstract}



\maketitle


\section{\bf Introduction}
\par A pseudo-Riemannian metric $g$ on a smooth manifold $M$ is called a Ricci soliton if there exists a smooth vector field $V$, such that 
\begin{equation}\label{e:1.1a}
\frac{1}{2}\pounds_{V}g + Ric = \lambda g, 
\end{equation}
where $\pounds_{V}$, $Ric$ and $\lambda$ denotes the Lie derivative in the direction of $V$,  the Ricci tensor and a real number respectively. A Ricci soliton $g$ is said to be a shrinking, steady or expanding according to whether $\lambda>0$, $\lambda=0$ or $\lambda<0$, respectively. In the last years, the interest in studying Ricci solitons has considerably increased among theoretical physicists in relation with string theory, and the fact that equation (\ref{e:1.1a}) is a special case of the Einstein field equations. Ricci solitons were introduced in Riemannian geometry \cite{Hamilton} as the self-similar solutions of the Ricci flow, and play an important role in understanding its singularities. A wide survey on Riemannian Ricci solitons may be found in \cite{Cao}. After their introduction in Riemannian case, the study of pseudo-Riemannian Ricci solitons attracted a growing number of authors. 
\par In parallel with contact and complex structures in the Riemannian case, paracontact metric structures were introduced by S. Kaneyuki and F. L. Williams  \cite{kk} in pseudo-Riemannian settings, as a natural odd-dimensional counterpart to para Hermitian structures.  A systematic study of paracontact metric manifolds started with the paper \cite{zamkov}. The technical apparatus introduced in  \cite{zamkov} is essential for further investigations of paracontact metric geometry. The problem of studying Ricci solitons in the context of paracontact metric geometry was initiated by G. Calvaruso and D. Perrone \cite{Calvaruso}. The case of Ricci solitons in three - dimensional paracontact  geometry was treated by C. L. Bejan and M. Crasmareanu in  \cite{Crasmareanu2} respectively G. Calvaruso and A. Perrone in \cite{CP}. For some recent results and further references on pseudo-Riemannian Ricci solitons, we may refer to \cite{Blaga3, BCGG, CF, CZ ,PT} and references therein.
\par Recently, G. Kaimakamis and K. Panagiotidou \cite{KP} initiated the notion of $*$-Ricci soliton where they essentially modified the definition of Ricci soliton by replacing the Ricci tensor $Ric$ in (\ref{e:1.1a}) with the $*$-Ricci tensor $Ric^*$. A pseudo-Riemannian metric $g$ on a smooth manifold $M$ is called a $*$-Ricci soliton if there exists a smooth vector field $V$, such that 
\begin{equation}\label{e:1.1}
\frac{1}{2}(\pounds_{V}g)(X, Y) + Ric^{*}(X, Y) = \lambda g(X, Y), 
\end{equation}
where 
\begin{equation}\label{e:1.2}
Ric^{*}(X,Y)=\frac{1}{2}(trace\{\phi\cdot R(X, \phi Y)\}), 
\end{equation}
for all vector fields $X, Y$ on $M$. Here, it is mentioned that the notion of $*$-Ricci tensor was first introduced by S. Tachibana \cite{Tachibana} on almost Hermitian manifolds and further studied by T. Hamada \cite{Hamada} on real hypersurfaces of non-flat complex space forms. 
\par Motivated by the above mentioned works, in this frame-work, we make an effort to study a $*$ -Ricci soliton on paracontact geometry mainly concerned the special case of para-Sasakian manifold. Para-Sasakain manifolds have been studied in recent years by many authors, emphasizing similarities and differences with respect to the most  well known Sasakian case. In this junction, it is suitable to mention that, Ricci silitons on para-Sasakian manifolds were studied in the paper \cite{DGP}. In the present paper, our main object is to study $*$-Ricci soliton within the frame-work of para-Sasakian manifold and prove the following result. 
\begin{theorem}\label{thm1}
Let $M(\phi, \xi,\eta, g)$ be a $(2n+1)$-dimensional para-Sasakian manifold. If $g$ is a $*$-Ricci soliton on $M$, then either $M$ is $\mathbb{D}$-homothetic to an Einstein manifold, or the Ricci tensor of $M$ with respect to canonical paracontact connection vanishes. In the first case, the soliton vector field is Killing and in the second case, the soliton vector field leaves $\phi$ invariant.
\end{theorem}
\section{\bf Preliminaries}
In this section we collecting some basic definitions and formulas on paracontact metric manifolds and para-Sasakian manifolds. All the manifolds are assumed to be connected and smooth. We may refer to \cite{ACGL, CKM, SS, DGP, VD, zamkov} and references therein for more information about para-Sasakian geometry.
\par An almost paracontact structure on a $(2n + 1)$-dimensional (connected) smooth manifold $M$ is a triple $(\phi, \xi, \eta)$, where $\phi$ is a $(1, 1)$-tensor, $\xi$ a global vector field and $\eta$ a 1-form, such that 
\begin{equation}\label{e:2.1}
 \phi(\xi)=0, \,\,\,  \eta \cdot \phi=0, \,\,\,  \eta(\xi)=1, \,\,\, \phi^{2}=Id-\eta\otimes\xi(1),
\end{equation}
and the restriction $J$ of $\phi$ on the horizontal distribution ker$\eta$ is an almost paracomplex structure (that is, the eigensubbundles $D^+, D^-$ corresponding to the eigenvalues $1$,$-1$ of $J$ have equal dimension $n$). A pseudo-Riemannian metric $g$ on $M$ is {\it compatible} with the almost paracontact structure $(\phi, \xi, \eta)$ when
\begin{equation}\label{e:2.2}
g(\phi X, \phi Y) = −g(X, Y) + \eta(X)\eta(Y)
\end{equation}
In such a case, $(\phi, \xi, \eta, g)$ is said to be an almost paracontact metric structure. Remark that, by (\ref{e:2.1}) and (\ref{e:2.2}), $\eta(X) = g(\xi, X)$ for any compatible metric. Any almost paracontact structure admits compatible metrics, which, by (\ref{e:2.2}), necessarily have signature $(n+1, n)$. The {\it fundamental 2-form} $\phi$ of an almost paracontact metric structure $(\phi, \xi, \eta, g)$ is defined by $\phi(X, Y)=g(X, \phi Y)$, for all tangent vector fields $X, Y$. If $\phi=d\eta$, then the manifold $(M, \eta, g)$ (or $M(\phi, \xi, \eta, g)$) is called a {\it paracontact metric manifold and g the associated metric}. If the paracontact metric structure $M(\phi, \xi, \eta, g)$ is normal, that is, satisfies $[\phi, \phi]+2d\eta\otimes\xi=0$, then $(\phi, \xi, \eta, g)$ is  called para-Sasakian. Equivalently, a para contact metric structure $(\phi, \xi, \eta, g)$ is para-Sasakian if
\begin{equation}\label{e:2.3}
(\nabla_X \phi)Y=-g(X, Y)\xi+\eta(Y)X,
\end{equation}
for any vector fields $X$, $Y$ on $M$, where $\nabla$ is Levi-Civita connection of $g$. Alternatively, a paracontact metric structure on $M$ is said to be para-Sasakian if the metric cone is para-Kaehler \cite{ACGL}. Any para-Sasakian manifold is $K$-paracontact, and the converse also holds when $n=1$, that is for three-dimensional space. From (\ref{e:2.3}), it follows that
\begin{equation}\label{e:2.4}
\nabla_X \xi=-\phi X.
\end{equation}
Also in a $(2n+1)$-dimensional para-Sasakian manifold, the following  relations hold:
\begin{eqnarray}
\label{e:2.5} R(X,Y)\xi&=&\eta(X)Y-\eta(Y)X,\\
\label{e:2.6} R(X,\xi)\xi&=&-X+\eta(X)\xi,\\
\label{e:2.7} Ric(X,\xi)&=&-2n\eta(X),\\
\label{e:2.8} Q\xi&=&-2n\xi,
\end{eqnarray}
for any $X, Y$ on $M$. Here, $R$ denotes the curvature tensor of $g$ and $Ric$ denotes the Ricci tensor defined by $S(X,Y)=g(QX, Y)$, where $Q$ is the Ricci operator. 
\begin{lemma}\label{l:2.1}
Let $M(\phi, \xi, \eta, g)$ be a para-Sasakian manifold. Then
\begin{center}
(i) $\nabla_{\xi}Q=0$,  and (ii) $(\nabla_{X} Q)\xi=Q\phi X+2n\phi X$.
\end{center} 
\end{lemma} 
{\bf Proof:} Since $\xi$ is Killing, we have $\pounds_V Ric=0$. This implies ($\pounds_{\xi} Q)X=0$ for any vector field $X$ on $M$. From which it follows that 
\begin{eqnarray}\nonumber
0&=&\pounds_{\xi}(QX)-Q(\pounds_\xi X)\\
&=&\nabla_{\xi} QX+\nabla_{QX}\xi-Q(\nabla_\xi X)+Q(\nabla_X \xi)\nonumber\\
&=&(\nabla_\xi Q)X+\nabla_{QX} \xi +Q(\nabla_X\xi)\nonumber.
\end{eqnarray}
Using (\ref{e:2.4}) in the above equation gives $\nabla_{\xi}Q=Q\phi-\phi Q$.
Since the Ricci operator $Q$ commutes with $\phi$ on para-Sasakian manifold, we have $(i)$.
Next, taking covariant differentiation of (\ref{e:2.8}) along an arbitrary vector field $X$ on $M$ and using (\ref{e:2.4}), we obtain $(ii)$. This completes the proof.
\par If the Ricci tensor of a para-Sasakian manifold $M$ is of the form
\begin{equation}
\nonumber Ric(X,Y)=A g(X,Y)+B\eta(X)\eta(Y),
\end{equation}
 for any vector fields $X, Y$ on $M$, where $A$ and $B$ being constants, then $M$ is called an $\eta$-Einstein manifold.
\par The 1-form $\eta$ is determined up to a horizontal distribution and hence $\mathbb{D}=Ker\eta$ are connected by $\widetilde{\eta}=\sigma \eta$ for a positive smooth function $\sigma$ on a paracontact manifold $M$. This paracontact form $\bar{\eta}$ defines the structure tensor $(\bar{\phi}, \bar{\xi}, \bar{g})$ corresponding to $\widetilde{\eta}$ using the condition given in the paper \cite{zamkov}. We call the transformation of the structure tensors given by Lemma 4.1 of \cite{zamkov} a gauge (conformal) transformation of paracontact pseudo-Riemannian structure. When $\sigma$ is constant this is a $\mathbb{D}$-homothetic transformation. Let $M(\phi, \xi, \eta, g)$ be a paracontact manifold and 
\begin{equation}
\nonumber \bar{\phi}=\phi,\,\,\, \bar{\xi}=\frac{1}{\alpha}\xi,\,\,\,\bar{\eta}=\alpha\eta,\,\,\, \bar{g}=\alpha g+(\alpha^2-\alpha)\eta\otimes\eta,  \,\,\, \alpha=const.\neq 0
\end{equation}
 to be $\mathbb{D}$-homothetic transformation. Then $(\bar{\phi}, \bar{\xi}, \bar{\eta}, \bar{g})$ is also a paracontact structure. Here the para-Sasakian structure is preserved since the normality conditions is preserved under $\mathbb{D}$-homothetic transformations. Using the formula appeared in \cite{zamkov} for $\mathbb{D}$-homothetic deformation, one can easily verify that if $M(\phi, \xi, \eta, g)$ is a (2n+1) - dimensional $(n>1)$ $\eta$-Einstein para-Sasakian structure with scalar curvature $r \neq 2n$, then there exists a constant $\alpha$ such that $M(\bar{\phi}, \bar{\xi}, \bar{\eta}, \bar{g})$  is an Einstein para-Sasakian structure. So we adopt the following result.
 \begin{lemma}
 Any $(2n+1)$-dimensional $\eta$-Einstein para-Sasakian manifold with scalar curvature not equal to $2n$ is $\mathbb{D}$-homothetic to an Einstein manifold.
\end{lemma}
\par The canonical paracontact connection on a paracontact manifold was defined by S. Zamkovoy \cite{zamkov}. On an integrable paracontact metric manifold such  a connection is unique and is defined in terms of the Levi-Civita connection by
\begin{equation}\label{e:2.9}
\widetilde{\nabla}_X Y=\nabla _X Y+\eta (X)\phi Y-\eta(Y)\nabla _X \xi+(\nabla_X \eta)(Y)\xi.
\end{equation}
\par A canonical paracontact connection on a para-Sasakian manifold which seems to be the paracontact analogue of the (generalized) Tanaka-Webster connection.
\par In view of (\ref{e:2.4}) in (\ref{e:2.9}), we arrive at
 \begin{equation}\label{e:2.10}
\widetilde{\nabla}_X Y=\nabla _X Y+\eta (X)\phi Y+\eta(Y)\phi X+g(X, \phi Y)\xi
\end{equation}
for all vector fields $X$, $Y$ on $M$. The connection $\widetilde{\nabla}$ given by (\ref{e:2.10}) is called a canonical paracontact connection on a para-Sasakian manifold. Indeed, the Ricci tensor $\widetilde{Ric}$ of a $(2n+1)$-dimensional para-Sasakian manifold with respect to canonical paracontact connection $\widetilde{\nabla}$ is defined by
\begin{equation}\label{e:2.11}
\widetilde{Ric} (X,Y)=Ric(X,Y)-2g(X,Y)+(2n+2)\eta(X)\eta(Y)
\end{equation}
for all vector fields $X$, $Y$ on $M$. For details we refer to \cite{IVZ} and \cite{zamkov}. 
\section{\bf Proof of Theorem \ref{thm1}}
\par In this section, before presenting our main result about $*$-Ricci soliton on a para-Sasakian manifold, we state and prove some lemmas which will be used to prove Theorem \ref{thm1}.
\begin{lemma}\label{l:3.1}
The $*$-Ricci tensor on a $(2n+1)$-dimensional para-Sasakian manifold $M(\phi, \xi, \eta, g)$ is given by
\begin{equation}\label{e:3.1}
Ric^{*}(X,Y)=-Ric(X,Y)-(2n-1)g(X,Y)-\eta(X)\eta(Y)
\end{equation}
for any vector fields $X, Y$ on $M$.
\end{lemma}
{\bf Proof:} The Ricci tensor $Ric$ of a $(2n+1)$-dimensional para-Sasakian manifold $M(\phi, \xi, \eta, g)$  satisfies the relation (c.f. Lemma 3.15 in \cite{zamkov}):
\begin{eqnarray}\label{e:3.2}
Ric(X,Y)&=&\frac{1}{2}\sum_{i=1}^{2n+1}R'(X, \phi Y, e_i, \phi e_i)\nonumber\\
&-&(2n-1)g(X,Y)-\eta(X)\eta(Y)
\end{eqnarray}
for any vector fields $X, Y$ on $M$. Using the skew-symmetric property of $\phi$, we write
\begin{equation}\nonumber
\sum_{i=1}^{2n+1}R'(X, \phi Y, e_i, \phi e_i)=\sum_{i=1}^{2n+1}(R(X, \phi Y), e_i, \phi e_i)=\sum_{i=1}^{2n+1}g(\phi R(X, \phi Y), e_i, e_i).
\end{equation}
 By virtue of this, it follows from (\ref{e:3.2}) that
 \begin{equation}\label{e:3.3}
  \sum_{i=1}^{2n+1}g(\phi R(X, \phi Y) e_i, e_i)=-2 Ric(X,Y)-2(2n-1)g(X,Y)-2\eta(X)\eta(Y).
 \end{equation}
 Making use of (\ref{e:1.2}) in (\ref{e:3.3}), we obtain (\ref{e:3.1}).
\begin{lemma}\label{l:3.2}
For a para-Sasakian manifold, we have the following relation
\begin{equation}\label{e:3.4}
(\pounds_V \eta)(\xi)=-\eta(\pounds_V \xi)=\lambda.
\end{equation}
\end{lemma}
{\bf Proof:} By virtue of Lemma \ref{l:3.1}, the $*$-Ricci soliton equation (\ref{e:1.1}) can be expressed as 
\begin{equation}\label{e:3.5}
(\pounds_{V}g)(X,Y)=2Ric(X,Y)+2(2n-1+\lambda)g(X,Y)+2\eta(X)\eta(Y).
\end{equation}
Setting $Y=\xi$ in (\ref{e:3.5}) and using (\ref{e:2.7}) it follows that $(\pounds_{V}g)(X,\xi)=2\lambda \eta(X)$. Lie-differentiating the equation $\eta(X)=g(X, \xi)$ along $V$ and by virtue of last equation, we find
\begin{equation}\label{e:3.6}
(\pounds_{V}\eta)(X)-g(\pounds_{V}\xi, X)-2\lambda \eta(X)=0.
\end{equation}
Next, Lie-derivative of $g(\xi, \xi)=1$ along $V$ and equation (\ref{e:3.6}) completes proof.
\begin{lemma}\label{l:3.3}
Let $M(\phi, \xi, \eta, g)$ be a $(2n+1)$-dimensional para-Sasakian manifold. If $g$ is a $*$-Ricci soliton, then $M$ is an $\eta$-Einstein manifold and the Ricci tensor can be expressed as 
\begin{equation}\label{e:3.17}
Ric(X,Y)=-\left[2n-1+\frac{\lambda}{2}\right] g(X,Y)+\left[\frac{\lambda}{2}-1\right]\eta(X)\eta(Y).
\end{equation}
for any vector fields $X, Y$ on $M$.
\end{lemma}
{\bf Proof:} First, taking covariant differentiation of (\ref{e:3.5}) along an arbitrary vector field $Z$, we get
\begin{eqnarray}\label{e:3.7}
&&(\nabla_Z \pounds_{V}g)(X,Y)\nonumber\\
&=& 2\{(\nabla_Z Ric)(X,Y)-g(X, \phi Z)\eta(Y)-g(Y, \phi Z)\eta(X)\}.
\end{eqnarray}
According to Yano \cite{KY}, we get
\begin{eqnarray}\nonumber
&&(\pounds_{V}  \nabla _Z g - \nabla _Z \pounds_{V} g - \nabla_{[V, Z]g})(X,Y)\\
&=& -g((\pounds_V \nabla)(Z,X),Y)-g((\pounds_V \nabla)(Z,Y),X), \nonumber
\end{eqnarray}
for any vector fields $X, Y, Z$ on $M$. In view of the parallelism of the pseudo-Riemannian metric $g$, we get from the above relation that
\begin{equation}\label{e:3.8}
(\nabla_Z \pounds_{V}g)(X,Y)=g((\pounds_{V} \nabla)(Z, X),Y)+g((\pounds_{V} \nabla)(Z, Y),X).
\end{equation}
Comparing (\ref{e:3.7}) and (\ref{e:3.8}), we have
\begin{eqnarray}
\label{e:3.9}
&&g((\pounds _V \nabla)(Z,X),Y)+g((\pounds _V \nabla)(Z,Y),X)\nonumber\\
&=&2\{(\nabla_Z Ric)(X,Y)-g(X,\phi Z)\eta(Y)-g(Y, \phi Z)\eta(X)\}.
\end{eqnarray}
By a straightforward combinatorial combination of (\ref{e:3.9}) gives 
\begin{eqnarray}
\label{e:3.10}
g((\pounds_V \nabla)(X,Y),Z)&=&-(\nabla_Z Ric)(X,Y)+(\nabla_X Ric)(Y,Z)\nonumber\\
&+&(\nabla_Y Ric)(Z,X) + 2g(X, \phi Z)\eta(Y)\nonumber\\
&+& 2g(Y, \phi Z)\eta(X).
\end{eqnarray}
Replacing $Y$ by $\xi$ in (\ref{e:3.10}) and Lemma \ref{l:2.1}, we have
\begin{equation}\label{e:3.11}
(\pounds_{V} \nabla)(X,Y)=2(2n-1)\phi X+2Q\phi X.
\end{equation}
Further, differentiating (\ref{e:3.11}) covariantly along an arbitrary vector field $Y$ on $M$ and then using the relations (\ref{e:2.3}) and (\ref{e:2.8}), we get
\begin{eqnarray}\label{e:3.12}
&&(\nabla_Y\pounds_{V}\nabla)(X, \xi)+(\pounds_{V}\nabla)(X, \phi Y)\\
&=&2\{(\nabla_Y Q)\phi X+\eta(X)QY+(2n-1)\eta(X)Y+g(X,Y)\xi\}.\nonumber
\end{eqnarray}
Again, according to Yano \cite{KY} we have the following commutation formula  
\begin{equation}\label{e:3.13}
(\pounds_{V}R)(X,Y)Z=(\nabla_X \pounds_{V}\nabla)(Y,Z)-(\nabla_Y \pounds_{V}\nabla)(X,Z).
\end{equation}
Replace $Z$ by $\xi$ in (\ref{e:3.13}) and taking into account of (\ref{e:3.12}), we obtain
\begin{eqnarray}\label{e:3.14}
&&(\pounds_{V}R)(X, Y)\xi+(\pounds_{V}\nabla)(Y, \phi X)-(\pounds_{V}\nabla)(X, \phi Y)\nonumber\\
&=&2\{(\nabla_X Q)\phi Y-(\nabla_Y Q)\phi X+\eta(Y)QX-\eta(X)QY\nonumber\\
&+&(2n-1)(\eta(Y)X-\eta(X)Y)\}.
\end{eqnarray}
Taking $\xi$ for $Y$ in (\ref{e:3.14}), then using (\ref{e:2.8}), (\ref{e:3.11}) and Lemma      
 \ref{l:2.1}, we have
\begin{equation}\label{e:3.15}
(\pounds_{V}R)(X, \xi)\xi=4\{QX+(2n-1)X+\eta(X)\xi\}.
\end{equation} 
Next, taking Lie-derivative of (\ref{e:2.6}) along $V$ and taking into account of (\ref{e:2.5}) and (\ref{e:3.4}) one can get
\begin{equation}\label{e:3.16}
(\pounds_{V}R)(X, \xi)\xi=(\pounds_{V} \eta)(X)\xi-g(\pounds_{V} X, \xi)-2\lambda X.
\end{equation} 
Comparing (\ref{e:3.15}) with (\ref{e:3.16}), and making use of (\ref{e:3.6}), we obtain the required result.\\
\par Now, we state and prove the main result in the following:\\
\
\\
{\bf Proof of Theorem \ref{thm1}:} Making use of (\ref{e:3.17}), the soliton equation (\ref{e:3.5}) takes the form
\begin{equation}\label{e:3.17b}
(\pounds_{V} g)(X,Y)=\lambda\{g(X,Y)+\eta(X)\eta(Y)\}.
\end{equation}
Taking Lie-differentiation of (\ref{e:3.17}) along the vector field $V$ and using (\ref{e:3.5}) we obtain
\begin{eqnarray}\label{e:3.17a}
(\pounds_{V} Ric)(X,Y)&=&\left(\frac{\lambda}{2}-1\right)\{\eta(Y)(\pounds_{V}\eta)(X)+\eta(X)(\pounds_{V}\eta)(Y)\}\nonumber\\
&-&\left(2n-1+\frac{\lambda}{2}\right)\lambda\{g(X,Y)+\eta(X)\eta(Y)\}. 
\end{eqnarray}
Other hand, differentiating the equation (\ref{e:3.17}) covariantly along an arbitrary vector field $Z$ on $M$ and then using (\ref{e:2.4}) we find
\begin{equation}\label{e:3.18} 
(\nabla_{Z} Ric)(X,Y)=\left(1-\frac{\lambda}{2}\right)\{g(X, \phi Z)\eta(Y)+g(Y, \phi Z)\eta(X)\}.
\end{equation}
In view of (\ref{e:3.18}), equation (\ref{e:3.10}) transforms into
\begin{equation}\label{e:3.19}
(\pounds_{V} \nabla)(X,Y)=-\lambda\{\eta(Y)\phi X+\eta(X)\phi Y\}.
\end{equation}
Differentiating (\ref{e:3.19}) covariantly along an arbitrary vector field $Z$ on $M$ and making use of (\ref{e:2.3}) and (\ref{e:2.4}) yields
\begin{eqnarray}\label{e:3.20}
\nonumber (\nabla_Z\pounds_{V} \nabla)(X,Y)&=&\lambda\{g(Y, \phi Z)\phi X+g(X, \phi Z)\phi Y+g(X,Z)\eta(Y)\xi\nonumber\\
&+&g(Y,Z)\eta(X)\xi-2\eta(X)\eta(Y)Z \}.
\end{eqnarray}
Using (\ref{e:3.20}) in commutation formula (\ref{e:3.13}) and using (\ref{e:2.4}) we produce
\begin{eqnarray}\label{e:3.21}
\nonumber (\pounds_{V} R)(X,Y)Z&=&\lambda\{g(\phi X, Z)\phi Y- g(\phi Y, Z)\phi X+2g(\phi X, Y)\phi Z \nonumber\\
&+&g(X,Z)\eta(Y)\xi- g(Y,Z)\eta(X)\xi-2\eta(Y)\eta(Z)X\nonumber\\
&+&2\eta(X)\eta(Z)Y\}.
\end{eqnarray}
Contracting (\ref{e:3.21}) over $Z$, we have
\begin{equation}\label{e:3.22}
(\pounds_{V} Ric)(Y,Z)=2\lambda\{g(Y,Z)-(2n+1)\eta(Y)\eta(Z)\}.
\end{equation}
Comparison of (\ref{e:3.17a}) and (\ref{e:3.22}) gives
\begin{eqnarray}\label{e:3.23}
&&\left(\frac{\lambda}{2}-1\right)\{\eta(Y)(\pounds_{V} \eta)(Z)+\eta(Z)(\pounds_{V} \eta)(Y)\}	\nonumber\\
&-&\left(2n-1+\frac{\lambda}{2}\right)\lambda\{g(Y,Z)+\eta(Y)\eta(Z)\}\nonumber\\
&=&2\lambda\{g(Y,Z)-(2n+1)\eta(Y)\eta(Z)\}.
\end{eqnarray}
Replacing $Y$ by $\phi^2 Y$ in (\ref{e:3.23}) and then using (\ref{e:2.1}) and (\ref{e:3.4}) we have
\begin{equation}\label{e:3.24}
\left(\frac{\lambda}{2}-1\right)(\pounds_{V}\eta)(Y)\eta(Z)=\lambda\left[1+2n+\frac{\lambda}{2}\right]g(Y,Z)-2n\lambda\, \eta(Y)\eta(Z).
\end{equation}
Taking account of (\ref{e:3.24}) in (\ref{e:3.23}) and then plugging $Z$ by $\phi Z$, we obtain 
\begin{equation}\label{e:3.25}
\lambda\left[2n+1+\frac{\lambda}{2}\right]g(Y,\phi Z)=0.
\end{equation}
Since $\phi(Y, Z)=g(Y, \phi Z)$ is non-vanishing everywhere on $M$, we can obtain either $\lambda=0$ or $\lambda=-2(2n+1)$.\\
{\bf Case I:} 
If $\lambda=0$, from (\ref{e:3.17b}) we see that $\pounds _V g=0$, that is, $V$ is Killing. From (\ref{e:3.17}) we see that
\begin{equation}\label{e:3.26}
Ric(X, Y)=-(2n-1)g(X,Y)-\eta(X)\eta(Y).
\end{equation}
Contracting the equation (\ref{e:3.26}) we get $r=-4n^2$, where $r$ is the scalar curvature of the manifold $M$. This shows that $M$ is an $\eta$-Einstein manifold with scalar curvature $r\neq 2n$. Hence, we conclude that $M$ is $\mathbb{D}$-homothetic to an Einstein manifold.\\
{\bf Case II:} If $\lambda=-2(2n+1)$, then replacing $Z$ by $\xi$ in (\ref{e:3.24}) and then setting $Y$ by $\phi Y$ of the resulting equation gives 
$(\frac{\lambda}{2}-1)(\pounds _V \eta)(\phi Y)=0$. Since $\lambda=-2(2n+1)$, we have $\lambda\neq 2$. thus we have $(\pounds _V \eta)(\phi Y)=0$. Taking $Y$ by $\phi Y$ in the foregoing equation and using (\ref{e:2.1}), we obtain
\begin{equation}\label{e:3.28}
(\pounds _V \eta)(Y)=-2(2n+1)\eta(X).
\end{equation}
Taking exterior differentiation $d$ on (\ref{e:3.28}) we have 
\begin{equation}\label{e:3.29}
(\pounds _V d\eta)(X,Y)=-2(2n+1)g(X, \phi Y),
\end{equation}
noting that $d$ commutes with $\pounds_V$. Further, taking the Lie-derivative of the well known equation $d\eta(X,Y)=g(X, \phi Y)$ along the soliton vector field $V$ provides
\begin{equation}\label{e:3.30}
(\pounds _V d\eta)(X,Y)=(\pounds _V g)(X,\phi Y)+g(X, (\pounds _V \phi)Y).
\end{equation}
From (\ref{e:3.17b}) we also deduce
\begin{equation}\label{e:3.31}
(\pounds _V g)(X,\phi Y)=-2(2n+1)g(X, \phi).
\end{equation}
Using (\ref{e:3.29}) and (\ref{e:3.31}) in (\ref{e:3.30}) we find $\pounds _V \phi=0$. Hence, the soliton vector field $V$ leaves $\phi$ invariant. \\
Further, using $\lambda=-2(2n+1)$ in (\ref{e:3.17}) it follows that 
\begin{equation}\label{e:3.27}
Ric(X, Y)=2g(X,Y)-(2n+2)\eta(X)\eta(Y).
\end{equation}
Contracting (\ref{e:3.27}) we obtain $r=2n$ (i.e., the manifold $M$ cannot be $\mathbb{D}$-homothetic to an Einstein manifold). Then, by taking account of (\ref{e:3.27}) in (\ref{e:2.11}) we obtain $\widetilde{Ric}(X,Y)=0$. That is, the Ricci tensor with respect to the connection $\widetilde{\nabla}$ vanishes. This completes the proof of our theorem \ref{thm1}.
\section{\bf Conclusion}
The study of Ricci solitons on Riemannian manifolds and pseudo-Riemannian manifolds is an issue, which is of great importance in the area of differential geometry and in physics as well. Ricci soliton generalizes the notion of Einstein metric on a Riemannian manifold. The $*$-Ricci soliton is a new notion not only in the area of differential geometry but in the area of physics as well. Note that, a $*$-Ricci soliton is trivial if the vector field $V$ is Killing, and in this case the manifold becomes $*$-Einstein, i.e., $Ric^* = \lambda g$. So far we know that, the notion of $*$-Ricci tensor appears on complex and contact manifolds only. However, some classifications are available in the literature in terms of the $*$-Ricci tensor. On the other hand, in the recent years, many authors have studied and pointed out the importance of paracontact geometry and, in particular, of para-Sasakian geometry, by giving the relationships with the theory of para-Kahler manifolds and its role in pseudo-Riemannian geometry and mathematical physics. Here, making use of the formulas of para-Sasakian manifold one can easily deduce an expression of the $*$-Ricci tensor. Due to the pressence of some extra terms in the expression of $*$-Ricci tensor the defining condition of the $*$-Ricci soliton is different from Ricci soliton. Thus, in this connection, we are interested and studied $*$-Ricci solitons within the framework of para-Sasakian manifold. The results obtained in this paper are playing an important role in differential geometry and mathematical physics. 
\\
{\bf Acknowledgment:} The first author (DGP) is thankful to University Grants Commission, New Delhi, India, for financial support to the Department of Mathematics, K. U. Dharwad in the form of UGC-SAP-DRS-III programme.

\end{document}